\documentclass[draft,12pt]{amsart}
\usepackage{amssymb,times}
\usepackage[numbers]{natbib}

\newtheorem*{Thm*}{Theorem}
\newtheorem{Thm}{Theorem}

\newtheorem{Prop}[Thm]{Proposition}
\newtheorem{Lemma}{Lemma}

\theoremstyle{definition}
\newtheorem{Defn}{Definition}
\newtheorem{Notation}[Defn]{Notation}

\newtheorem{Remark}{Remark}
\newtheorem{Ex}[Remark]{Example}

\newcommand{\mf}[1]{\mathbb{#1}}
\newcommand{\mc}[1]{\mathcal{#1}}
\newcommand{\mb}[1]{\mathbf{#1}}
\newcommand{\mk}[1]{\mathfrak{#1}}

\DeclareMathOperator{\Falg}{\mathcal{F}_{\mathrm{alg}}}

\newcommand{\norm}[1]{\left\Vert#1\right\Vert}
\newcommand{\abs}[1]{\left\vert#1\right\vert}

\newcommand{\set}[1]{\left\{#1\right\}}
\newcommand{\ip}[2]{\left \langle #1, #2 \right \rangle}
\newcommand{\state}[1]{\varphi \left[ #1 \right]}
\renewcommand{\phi}{\varphi}

\newcommand{\Span}[1]{\mathrm{Span} \left( #1 \right)}

\newcommand{\br}{\medskip\noindent}

\allowdisplaybreaks[1]

\title{Monic non-commutative orthogonal polynomials}
\author[M.~Anshelevich]{Michael Anshelevich}
\thanks{This work was supported in part by NSF grant DMS-0613195}
\address{Department of Mathematics, Texas A\&M University, College Station, TX 77843-3368}
\email{manshel@math.tamu.edu}
\date{\today}

\begin{document}

\begin{abstract}
Among all states on the algebra of non-commutative polynomials, we characterize the ones that have \emph{monic} orthogonal polynomials. The characterizations involve recursion relations, Hankel-type determinants, and a representation as a joint distribution of operators on a Fock space.
\end{abstract}

\maketitle

\section{Introduction}
\noindent
For a measure $\mu$ on $\mf{R}$ all of whose moments are finite, there are two standard ways to normalize the polynomials orthogonal with respect to $\mu$. One can take the polynomials to be monic,
\[
P_n(x) = x^n + \textsl{ lower order terms}.
\]
Or one can take them to be ortho\textsl{normal},
\[
\int_{\mf{R}} P_n(x)^2 \,d\mu(x) = 1.
\]
For a measure $\mu$ on $\mf{R}^d$, the situation is more subtle. One can always orthogonalize the subspaces of polynomials of different \emph{total} degree (so that one gets a family of pseudo-orthogonal polynomials). The most common approach is to work directly with these subspaces, without producing individual orthogonal polynomials; see, for example \cite{Dunkl-Xu}. One can also further orthogonalize the polynomials of the same total degree, for example to make them orthonormal \cite{Constantinescu,Banks-Const-OPS2}; however this requires a choice of an order on the monomials of the same degree, and there is no canonical choice of such order. The third approach, and (it is easy to see) the only one that will produce monic orthogonal polynomials, is to require that the pseudo-orthogonal polynomials obtained in the first step already be orthogonal. The price one pays is that this can be done only for some measures $\mu$.

\br
This paper is about orthogonal polynomials in \emph{non-commuting} variables. There is no difficulty with the definition. The usual orthogonal polynomials are obtained by starting with a measure $\mu$ on $\mf{R}^d$, thinking of $\mf{R}[x_1, x_2, \ldots, x_d]$ as a vector space with the (pre-)inner product
\[
\ip{P}{Q} = \int_{\mf{R}^d} P(\mb{x}) Q(\mb{x}) \,d\mu(\mb{x}),
\]
and applying the Gram-Schmidt procedure to the monomials $\set{x_{u(1)} x_{u(2)} \ldots x_{u(n)}}$. In the non-commutative case, one starts directly with a positive linear functional (state) $\phi$ on the algebra of non-commutative polynomials $\mf{R} \langle x_1, x_2, \ldots, x_n \rangle$, and orthogonalizes the monomials in non-commuting variables with respect to the inner product
\[
\ip{P}{Q} = \state{P^\ast(\mb{x}) Q(\mb{x})}.
\]
In particular, in \cite{AnsMulti-Sheffer}, I showed that monic polynomials in non-commuting variables orthogonal with respect to a faithful state satisfy a recursion relation, as orthogonal polynomials should. The question remained: which states have monic orthogonal polynomials? It is the first question answered in Theorem~\ref{Thm:Monic-states} of this paper, where moment conditions and Fock space representations of such states are provided. The second question, asked by the referee, was whether one needs the faithfulness condition. In this paper, that condition is removed using a new technique, namely a representation of the state as a joint distribution of some operators. Such a representation is closely related to a combinatorial way of representing moments as sums over lattice paths \cite{Flajolet}, but in the multivariate case I find the operator formulation more useful.

\br
An important class of states that have monic orthogonal polynomials are the free Meixner states, whose study I initiated in \cite{AnsMulti-Sheffer}. Applications of the techniques from the current paper to free Meixner states are considered in a companion paper \cite{AnsFree-Meixner}.

\section{Preliminaries}

\subsection{Polynomials}
\noindent
Let $\mb{x} = (x_1, x_2, \ldots, x_d)$ be a $d$-tuple of non-commuting variables. Let
\[
\mf{R}\langle \mb{x} \rangle = \mf{R}\langle x_1, x_2, \ldots, x_d \rangle
\]
be all the polynomials with real coefficients in these variables. \emph{Multi-indices} are elements
\[
\vec{u} \in \set{1, \ldots, d}^k
\]
for $k \geq 0$; for $\abs{\vec{u}} = 0$ denote $\vec{u}$ by $\emptyset$. Monomials in non-commuting variables $(x_1, \ldots, x_d)$ are indexed by such multi-indices:
\[
x_{\vec{u}} = x_{u(1)} \ldots x_{u(k)}.
\]
Note that our use of the term ``multi-index'' is different from the usual one, which is more suited for indexing monomials in commuting variables.

\br
For two multi-indices $\vec{u}, \vec{v}$, denote by $(\vec{u}, \vec{v})$ their concatenation. For $\vec{u}$ with $\abs{\vec{u}} = k$, denote
\[
(\vec{u})^{op} = (u(k), \ldots, u(2), u(1)).
\]
Define an involution on $\mf{R}\langle \mb{x} \rangle$ via the $\mf{R}$-linear extension of
\[
(x_{\vec{u}})^\ast = x_{(\vec{u})^{op}}.
\]

\br
A \emph{monic polynomial family} in $\mb{x}$ is a family $\set{P_{\vec{u}}(\mb{x})}$ indexed by all multi-indices
\[
\bigcup_{k=1}^\infty \set{\vec{u} \in \set{1, \ldots, d}^k}
\]
(with $P_{\emptyset} = 1$ being understood) such that\[
P_{\vec{u}}(\mb{x}) = x_{\vec{u}} + \textsl{lower-order terms}.
\]
Note that $P_{\vec{u}}^\ast \neq P_{(\vec{u})^{op}}$ in general.

\begin{Defn}
\label{Defn:State}
A \emph{state} on $\mf{R} \langle \mb{x} \rangle$ is a functional
\[
\phi: \mf{R} \langle x_1, x_2, \ldots, x_d \rangle \rightarrow \mf{R}
\]
that is linear, compatible with the $\ast$-operation, that is for any $P$,
\[
\state{P} = \state{P^\ast},
\]
unital, that is $\state{1} = 1$, and positive, that is for any $P$,
\[
\state{P^\ast P} \geq 0.
\]
A state is \emph{faithful} if in the preceding equation, the equality holds only for $P = 0$. Unless noted otherwise, the states in this paper are \emph{not} assumed to be faithful.

\br
The numbers $\state{x_{\vec{u}}}$ are called the \emph{moments} of $\phi$.
\end{Defn}

\begin{Remark}
\label{Remark:Cauchy-Schwartz}
A state $\phi$ induces the pre-inner product
\[
\ip{P}{Q}_\phi = \state{P^\ast Q} = \ip{Q}{P}_\phi
\]
and the seminorm
\[
\norm{P}_\phi = \sqrt{\state{P^\ast P}}.
\]
Throughout the paper, we will typically drop $\phi$ from the notation, and denote the inner product and norm it induces simply by $\ip{\cdot}{\cdot}$, $\norm{\cdot}$.

\br
Even if a state is not faithful, it is easy to check that the Cauchy-Schwartz inequality
\[
\abs{\ip{P}{Q}} \leq \norm{P} \norm{Q}
\]
still holds. In particular, if $\norm{P} = 0$, then for any $Q$, $\ip{P}{Q} = 0$, and if $\norm{P - P'} = 0$, then for any $Q$, $\ip{P}{Q} = \ip{P'}{Q}$.
\end{Remark}

\begin{Defn}
A state has a \emph{monic orthogonal polynomial system}, or MOPS, if for any multi-index $\vec{u}$, there is a monic polynomial $P_{\vec{u}}$ with leading term $x_{\vec{u}}$, such that these polynomials are orthogonal with respect to $\phi$, that is,
\[
\ip{P_{\vec{u}}}{P_{\vec{v}}} = 0
\]
for $\vec{u} \neq \vec{v}$.
\end{Defn}

\noindent
Note that the same abbreviation is used in~\cite{Dumitriu-MOPS} to denote a class of multivariate orthogonal polynomials systems, which is different from ours.

\begin{Lemma}
\label{Lemma:L2-equal}
Let $\set{Q_{\vec{u}}}$ and $\set{P_{\vec{u}}}$ be two monic polynomial families, and suppose that $\set{Q_{\vec{u}}}$ are orthogonal with respect to a state $\phi$. The following are equivalent:
\begin{enumerate}
\item
$\set{P_{\vec{u}}}$ are pseudo-orthogonal with respect to $\phi$, which means that
\[
\ip{P_{\vec{u}}}{P_{\vec{v}}} = 0
\]
whenever $\abs{\vec{u}} \neq \abs{\vec{v}}$.
\item
For each $\vec{u}$, $Q_{\vec{u}}$ and $P_{\vec{u}}$ are equal in $L^2(\phi)$, that is,
\[
\norm{Q_{\vec{u}} - P_{\vec{u}}} = 0.
\]
\end{enumerate}
If either of these statements holds, then in fact $\set{P_{\vec{u}}}$ are orthogonal.
\end{Lemma}

\begin{proof}
(a) $\Rightarrow$ (b).
Since $Q_{\vec{u}}$ and $P_{\vec{u}}$ are monic, $Q_{\vec{u}} - P_{\vec{u}}$ has degree less than $\abs{\vec{u}}$. Since both families are pseudo-orthogonal, $Q_{\vec{u}} - P_{\vec{u}}$ is orthogonal to all polynomials of degree less than $\abs{\vec{u}}$. Therefore it has norm zero.

\br
(b) $\Rightarrow$ (a).
Using Remark~\ref{Remark:Cauchy-Schwartz}, for $\vec{v} \neq \vec{u}$, $\ip{P_{\vec{u}}}{P_{\vec{v}}} = \ip{Q_{\vec{u}}}{P_{\vec{v}}} = \ip{Q_{\vec{u}}}{Q_{\vec{v}}} = 0$. Note that this in fact implies that $\set{P_{\vec{u}}}$ are orthogonal.
\end{proof}

\begin{Lemma}
\label{Lemma:Ideal}
Let $\phi$ be a state with MOPS $\set{P_{\vec{v}}}$. Suppose that for some $\vec{u}$, $\norm{P_{\vec{u}}} = 0$. Then for any $i$, $\norm{P_{(i, \vec{u})}} = 0$.
\end{Lemma}

\begin{proof}
Since $\set{P_{\vec{v}}}$ are monic, they form a basis for the vector space of polynomials. So we may write
\[
x_i P_{\vec{u}}(\mb{x}) = P_{(i, \vec{u})}(\mb{x}) + \sum_{\abs{\vec{v}} \leq \abs{\vec{u}}} \alpha_{\vec{v}} P_{\vec{v}}(\mb{x}).
\]
Then
\[
\begin{split}
\norm{P_{(i, \vec{u})}}^2
& = \state{P_{(i, \vec{u})}^\ast(\mb{x}) P_{(i, \vec{u})}(\mb{x})}
= \state{P_{(i, \vec{u})}^\ast(\mb{x}) x_i P_{\vec{u}}(\mb{x})} - \sum_{\abs{\vec{v}} \leq \abs{\vec{u}}} \alpha_{\vec{v}} \state{P_{(i, \vec{u})}^\ast(\mb{x}) P_{\vec{v}}(\mb{x})} \\
& = \ip{x_i P_{(i, \vec{u})}(\mb{x})}{P_{\vec{u}}(\mb{x})} - \sum_{\abs{\vec{v}} \leq \abs{\vec{u}}} \alpha_{\vec{v}} \ip{P_{(i, \vec{u})}(\mb{x})}{P_{\vec{v}}(\mb{x})}
= 0
\end{split}
\]
because of orthogonality and Remark~\ref{Remark:Cauchy-Schwartz}.
\end{proof}

\section{Monic orthogonal polynomials states}
\label{Section:MOPS}
\noindent
The goal of this paper is to prove Theorem~\ref{Thm:Monic-states}, which provides a number of equivalent conditions describing a class of states. Since these conditions come from quite different frameworks, we first describe two constructions. Among other things, these results now apply to not-necessarily-faithful states, thus answering a question of the referee of \cite{AnsMulti-Sheffer}.

\begin{Notation}
\label{Notation:Hankel}
Put an order on all multi-indices that is compatible with degree and otherwise arbitrary, for example the lexicographic one. Let $\abs{\vec{u}} = n$. Construct a $1 + d + d^2 + \ldots + d^{n-1} + 1 = \frac{d^n - 1}{d - 1} + 1$-dimensional square matrix $A_{\vec{u}}$, whose rows and columns are labeled by all the multi-indices of length less than $n$ and $\vec{u}$, and whose $(\vec{v},\vec{w})$ entry is $\ip{x_{\vec{v}}}{x_{\vec{w}}}$. Denote $h_{\vec{u}} = \det A_{\vec{u}}$, and $\mk{h}_n$ the determinant of the matrix $A_{\vec{u}}$ with the $\vec{u}$'th row and column removed (so that it depends on $n$ but not on $\vec{u}$). These are multivariate versions of Hankel determinants, (compare with \cite{Banks-Const-OPS2}). Note that $h_{\vec{u}}, \mk{h}_n$ do not depend on the chosen order on multi-indices. Let $M_{\vec{u}}(\mb{x})$ be the matrix $A_{\vec{u}}$ with the $\vec{u}$'th row changed, so that the $(\vec{u}, \vec{w})$ entry is $x_{\vec{w}}$. Finally, let $h_{\vec{v}, \vec{u}}$ be the determinant of $A_{\vec{u}}$ with the $\vec{u}$'th row changed, so that the $(\vec{u}, \vec{w})$ entry is $\ip{x_{\vec{v}}}{x_{\vec{w}}}$. In particular, $h_{\vec{u}, \vec{u}} = h_{\vec{u}}$.
\end{Notation}

\subsection{Fock space construction}
\label{Subsec:General-Fock}
Let $\mc{H} = \mf{C}^d$, with the canonical orthonormal basis $e_1, e_2, \ldots, e_d$. Define the (algebraic) full Fock space of $\mc{H}$ to be
\[
\Falg(\mc{H}) = \bigoplus_{k=0}^\infty \mc{H}^{\otimes k}
\]
Equivalently, $\Falg(\mc{H})$ is the vector space of non-commutative polynomials in $e_1, e_2, \ldots, e_d$. Following convention, we will denote the generating vector in $\mc{H}^{\otimes 0} = \mf{C}$ by $\Omega$ instead of $1$.

\br
For each $k \geq 1$ let $\mc{C}^{(k)}$ be an operator
\[
\mc{C}^{(k)}: \mc{H}^{\otimes k} \rightarrow \mc{H}^{\otimes k}.
\]
We identify each $\mc{C}^{(k)}$ with its $d^k \times d^k$ matrix in the basis $\set{e_{u(1)} \otimes \ldots e_{u(k)}}$. Assume that for each $k$, $\mc{C}^{(k)}$ is diagonal and $\mc{C}^{(k)} \geq 0$.

\br
For $i = 1, 2, \ldots, d$, define $a_i^+$ and $a_i^-$ to be the usual (left) free creation and annihilation operators,
\begin{align*}
a_i^+ & \left(e_{u(1)} \otimes e_{u(2)} \otimes \ldots \otimes e_{u(k)} \right) = e_i \otimes e_{u(1)} \otimes e_{u(2)} \otimes \ldots \otimes e_{u(k)}, \\
a_i^-      & (e_j) = \ip{e_i}{e_j} \Omega = \delta_{i j} \Omega, \\
a_i^-      & \left(e_{u(1)} \otimes e_{u(2)} \otimes \ldots \otimes e_{u(k)} \right) = \ip{e_i}{e_{u(1)}} e_{u(2)} \otimes \ldots \otimes e_{u(k)}.
\end{align*}
For each $i = 1, 2, \ldots, d$ and each $k \geq 0$, let $\mc{T}_i^{(k)}$ be an operator
\[
\mc{T}_i^{(k)}: \mc{H}^{\otimes k} \rightarrow \mc{H}^{\otimes k}.
\]
Assume that each $\mc{T}_i^{(k)}$ satisfies
\begin{equation}
\label{Transpose}
\left(\mc{T}_i^{(k)}\right)^t \mc{K}_{\mc{C}} = \mc{K}_{\mc{C}} \mc{T}_i^{(k)},
\end{equation}
where $\mc{K}_{\mc{C}}$ is the operator in equation~\eqref{Kernel} below, and $A^t$ is the transpose. We will denote by $\mc{T}_i$ and $\mc{C}$ the operators on $\Falg(\mc{H})$ acting as $\mc{T}_i^{(k)}$ and $\mc{C}^{(k)}$ on each component. Finally, denote
\[
\tilde{a}_i^- = a_i^- \mc{C}.
\]
Note that $a_i^- \Omega = \tilde{a}_i^- \Omega = 0$.

\br
On each $\mc{H}^{\otimes k}$ one has the usual inner product $\ip{\cdot}{\cdot}$ induced from $\mc{H}$. Define a new pre-inner product $\ip{\cdot}{\cdot}_{\mc{C}}$ by using the non-negative kernel
\begin{equation}
\label{Kernel}
\mc{K}_{\mc{C}} = (I \otimes I \ldots \otimes I \otimes \mc{C}^{(1)}) \ldots (I \otimes I \otimes \mc{C}^{(k-2)}) (I \otimes \mc{C}^{(k-1)}) \mc{C}^{(k)}
\end{equation}
where $I$ denotes the identity matrix of the appropriate size. That is,
\[
\ip{\zeta}{\eta}_{\mc{C}} = \ip{\zeta}{\mc{K}_{\mc{C}} \eta}.
\]
Put this pre-inner product on each component of $\Falg(\mc{H})$ and define the components to be orthogonal among themselves. Factor out the vectors in $\Falg(\mc{H})$ of norm zero, and denote the completion of the factor with respect to this inner product $\mc{F}_{\mc{C}}(\mc{H})$.

\begin{Lemma}
\label{Lemma:Degeneracy}
Let $S = \set{\eta \in \Falg(\mc{H}) | \norm{\eta}_{\mc{C}} = 0}$, so that $\mc{F}_{\mc{C}}(\mc{H})$ is the completion of $\Falg(\mc{H}) / S$.
\begin{enumerate}
\item
$S = \Span{e_{u(1)} \otimes \ldots e_{u(k)} | k > 0, \norm{e_{u(1)} \otimes \ldots e_{u(k)}}_{\mc{C}} = 0} = \ker \mc{K}_{\mc{C}}$.
\item
The operators $a_i^+, \mc{T}_i, \tilde{a}_i^-$ factor through to $\mc{F}_{\mc{C}}(\mc{H})$.
\end{enumerate}
\end{Lemma}

\begin{proof}
(a). Let $\eta \in S$,
\[
\eta = \sum_{k=0}^n \sum_{\abs{\vec{u}} = k} \eta_{\vec{u}} e_{u(1)} \otimes \ldots e_{u(k)}.
\]
Since each $\mc{C}^{(k)}$ is diagonal, we can denote $\mc{C} (e_{u(1)} \otimes \ldots e_{u(k)}) = C_{\vec{u}} e_{u(1)} \otimes \ldots e_{u(k)}$ and
\[
\mc{K}_{\mc{C}} (e_{u(1)} \otimes \ldots e_{u(k)}) = \mc{K}_{\mc{C}, \vec{u}} e_{u(1)} \otimes \ldots e_{u(k)}.
\]
Then
\[
\begin{split}
0 & = \norm{\eta}_{\mc{C}}
= \ip{\sum_{k=0}^n \sum_{\abs{\vec{u}} = k} \eta_{\vec{u}} e_{u(1)} \otimes \ldots e_{u(k)}}{\mc{K}_{\mc{C}} \sum_{j=0}^n \sum_{\abs{\vec{v}} = j} \eta_{\vec{v}} e_{v(1)} \otimes \ldots e_{v(j)}} \\
& = \ip{\sum_{k=0}^n \sum_{\abs{\vec{u}} = k} \eta_{\vec{u}} e_{u(1)} \otimes \ldots e_{u(k)}}{\sum_{j=0}^n \sum_{\abs{\vec{v}} = j} \mc{K}_{\mc{C}, \vec{v}} \eta_{\vec{v}} e_{v(1)} \otimes \ldots e_{v(j)}}
= \sum_{\vec{u}} \eta_{\vec{u}}^2 \mc{K}_{\mc{C}, \vec{u}}.
\end{split}
\]
Since $\mc{K}_{\mc{C}, \vec{u}} \geq 0$, it follows that
\[
\eta_{\vec{u}} \neq 0 \Rightarrow \mc{K}_{\mc{C}, \vec{u}} = 0 \Rightarrow \norm{e_{u(1)} \otimes \ldots e_{u(k)}}_{\mc{C}} = 0
\]
and also that $\mc{K}_{\mc{C}} (\eta) = 0$.

\br
(b). By part (a), it suffices to consider the actions on $e_{u(1)} \otimes \ldots e_{u(k)} \in S$. Then
\[
\begin{split}
\mc{K}_{\mc{C}} a_i^+ (e_{u(1)} \otimes \ldots e_{u(k)})
& = (I \otimes \mc{K}_{\mc{C}}) \mc{C} (e_i \otimes e_{u(1)} \otimes \ldots e_{u(k)}) \\
= C_{(i, \vec{u})} (I \otimes \mc{K}_{\mc{C}}) (e_i \otimes e_{u(1)} \otimes \ldots e_{u(k)})
& = C_{(i, \vec{u})} a_i^+ \mc{K}_{\mc{C}} (e_{u(1)} \otimes \ldots e_{u(k)})
= 0.
\end{split}
\]
Also $\mc{K}_{\mc{C}} \mc{T}_i = \mc{T}_i^t \mc{K}_{\mc{C}}$ by definition of $\mc{T}_i$, so that $\mc{T}_i (\ker \mc{K}_{\mc{C}}) \subset \ker \mc{K}_{\mc{C}}$. Finally, it is easy to check that
\begin{equation}
\label{a-K-commute}
\mc{K}_{\mc{C}} \tilde{a}_i^- = \mc{K}_{\mc{C}} a_i^- \mc{C} = a_i^- (I \otimes \mc{K}_{\mc{C}}) \mc{C} = a_i^- \mc{K}_{\mc{C}}
\end{equation}
and so $\tilde{a}^-_i(\ker \mc{K}_{\mc{C}}) \subset \ker \mc{K}_{\mc{C}}$.
\end{proof}

\br
For each $i$, define an operator $\mc{X}_i$ on $\mc{F}_{\mc{C}}(\mc{H})$ with dense domain $\Falg(\mc{H})/S$ by
\[
\mc{X}_i = a_i^+ + \mc{T}_i + \tilde{a}_i^-.
\]

\begin{Prop}
\label{Prop:Symmetric}
Each $\mc{X}_i$ is a symmetric operator.
\end{Prop}

\begin{proof}
Using equation~\eqref{a-K-commute}, for any $\vec{u}, \vec{v}$,
\[
\begin{split}
& \ip{a_i^+ \left(e_{u(1)} \otimes \ldots \otimes e_{u(k)} \right)}{e_{v(0)} \otimes \ldots \otimes e_{v(k)}}_{\mc{C}} \\
&\quad = \ip{e_i \otimes e_{u(1)} \otimes \ldots \otimes e_{u(k)}}{e_{v(0)} \otimes \ldots \otimes e_{v(k)}}_{\mc{C}} \\
&\quad = \ip{e_i \otimes e_{u(1)} \otimes \ldots \otimes e_{u(k)}}{\mc{K}_{\mc{C}} \left(e_{v(0)} \otimes \ldots \otimes e_{v(k)} \right)} \\
&\quad = \ip{e_{u(1)} \otimes \ldots \otimes e_{u(k)}}{a_i^- \mc{K}_{\mc{C}} \left(e_{v(0)} \otimes \ldots \otimes e_{v(k)} \right)} \\
&\quad = \ip{e_{u(1)} \otimes \ldots \otimes e_{u(k)}}{\mc{K}_{\mc{C}} \tilde{a}_i^- \left(e_{v(0)} \otimes \ldots \otimes e_{v(k)} \right)} \\
&\quad = \ip{e_{u(1)} \otimes \ldots \otimes e_{u(k)}}{\tilde{a}_i^- \left(e_{v(0)} \otimes \ldots \otimes e_{v(k)} \right)}_{\mc{C}},
\end{split}
\]
so with respect to the $\mc{C}$-inner product, $a_i^+ + \tilde{a}_i^-$ is symmetric. Similarly,
\[
\begin{split}
& \ip{\mc{T}_i \left(e_{u(1)} \otimes \ldots \otimes e_{u(k)} \right)}{e_{v(1)} \otimes \ldots \otimes e_{v(k)}}_{\mc{C}} \\
&\qquad = \ip{\mc{T}_i \left(e_{u(1)} \otimes \ldots \otimes e_{u(k)} \right)}{\mc{K}_{\mc{C}} \left(e_{v(1)} \otimes \ldots \otimes e_{v(k)} \right)} \\
&\qquad = \ip{e_{u(1)} \otimes \ldots \otimes e_{u(k)}}{\mc{T}_i^t \mc{K}_{\mc{C}} \left(e_{v(1)} \otimes \ldots \otimes e_{v(k)} \right)} \\
&\qquad = \ip{e_{u(1)} \otimes \ldots \otimes e_{u(k)}}{\mc{K}_{\mc{C}} \mc{T}_i \left(e_{v(1)} \otimes \ldots \otimes e_{v(k)} \right)} \\
&\qquad = \ip{e_{u(1)} \otimes \ldots \otimes e_{u(k)}}{\mc{T}_i \left(e_{v(1)} \otimes \ldots \otimes e_{v(k)} \right)}_{\mc{C}},
\end{split}
\]
so $\mc{T}_i$ is symmetric.
\end{proof}

\begin{Defn}
\label{Defn:Fock-state}
For any choice of the matrices $\mc{C}^{(k)}$ and $\mc{T}_i^{(k)}$ as above, the corresponding (Fock) state $\phi = \phi_{\mc{C}, \set{\mc{T}_i}}$ on $\mf{R} \langle \mb{x} \rangle$ is defined by
\begin{equation*}
\state{P(x_1, x_2, \ldots, x_d)} = \ip{\Omega}{P(\mc{X}_1, \mc{X}_2, \ldots, \mc{X}_d) \Omega}_{\mc{C}}.
\end{equation*}
\end{Defn}

\begin{Lemma}
$\phi$ satisfies all the properties in Definition~\ref{Defn:State}.
\end{Lemma}

\begin{proof}
Since $\phi$ is a vector state corresponding to the vector $\Omega$ and $\norm{\Omega} = 1$, $\phi$ is linear, unital, and positive. Finally, since each $\mc{X}_i$ is symmetric,
\[
\begin{split}
\state{x_{u(1)} \ldots x_{u(k)}}
& = \ip{\Omega}{\mc{X}_{u(1)} \ldots \mc{X}_{u(k)} \Omega}_{\mc{C}}
= \ip{\mc{X}_{u(k)} \ldots \mc{X}_{u(1)} \Omega}{\Omega}_{\mc{C}} \\
& = \ip{\Omega}{\mc{X}_{u(k)} \ldots \mc{X}_{u(1)} \Omega}_{\mc{C}}
= \state{(x_{u(1)} \ldots x_{u(k)})^\ast}.
\qedhere
\end{split}
\]
\end{proof}

\begin{Lemma}
\label{Lemma:Fock-polynomials}
For any $\mc{C}, \set{\mc{T}_i}$ as above, the state $\phi_{\mc{C}, \set{\mc{T}_i}}$ has a MOPS $\set{P_{\vec{u}}}$, which satisfy
\begin{equation}
\label{Wick-representation}
P_{\vec{u}}(\mc{X}_1, \mc{X}_2, \ldots, \mc{X}_d) \Omega = e_{u(1)} \otimes \ldots \otimes e_{u(\abs{\vec{u}})}.
\end{equation}
\end{Lemma}

\begin{proof}
For each $\vec{u}$,
\[
\mc{X}_{\vec{u}} \Omega = e_{u(1)} \otimes \ldots \otimes e_{u(\abs{\vec{u}})} + \eta_{\vec{u}}
\]
for some $\eta_{\vec{u}} \in \bigoplus_{i=0}^{\abs{\vec{u}} - 1} \mc{H}^{\otimes i}$. So we can recursively construct (non-unique) monic polynomials $P_{\vec{u}}(\mb{x})$ which satisfy equation~\eqref{Wick-representation}. Note that since all $\mc{X}_i$ are symmetric, $P_{\vec{u}}^\ast(\mb{X})$ really is the adjoint of $P_{\vec{u}}(\mb{X})$. Then
\[
\begin{split}
\ip{P_{\vec{u}}}{P_{\vec{v}}}
& = \phi_{\mc{C}, \set{\mc{T}_i}} \left[ P_{\vec{u}}^\ast(\mb{x}) P_{\vec{v}}(\mb{x}) \right]
= \ip{\Omega}{P_{\vec{u}}^\ast(\mb{X}) P_{\vec{v}}(\mb{X}) \Omega}_{\mc{C}} \\
& = \ip{P_{\vec{u}}(\mb{X}) \Omega}{P_{\vec{v}}(\mb{X}) \Omega}_{\mc{C}}
= \ip{e_{u(1)} \otimes \ldots \otimes e_{u(\abs{\vec{u}})}}{e_{v(1)} \otimes \ldots \otimes e_{v(\abs{\vec{v}})}}_{\mc{C}}
= 0
\end{split}
\]
for $\vec{u} \neq \vec{v}$. Thus $\set{P_{\vec{u}}}$ are monic orthogonal polynomials.
\end{proof}

\noindent
In the following theorem, the most interesting characterizations are part (g), which uses only moments of $\phi$, and part (e), which provides a way to construct any $\phi$ from matricial data.

\begin{Thm}
\label{Thm:Monic-states}
Let $\phi$ be a state on $\mf{R} \langle \mb{x} \rangle$. The following are equivalent:
\begin{enumerate}
\item
The state $\phi$ has a monic orthogonal polynomial system.
\item
The polynomials $\set{P_{\vec{u}}}$ defined recursively by the Gram-Schmidt relation
\begin{equation}
\label{Gram-Schmidt}
P_{\vec{u}} = x_{\vec{u}} - \sum_{\vec{v}:\ \abs{\vec{v}} < \abs{\vec{u}},\ \norm{P_{\vec{v}}} \neq 0} \frac{\ip{x_{\vec{u}}}{P_{\vec{v}}}}{\norm{P_{\vec{v}}}^2} P_{\vec{v}}
\end{equation}
are a monic orthogonal polynomial system for $\phi$.
\item
There is a family of polynomials $\set{P_{\vec{u}}}$ such that $\state{P_{\vec{u}}} = 0$ for all $\vec{u} \neq \emptyset$ and they satisfy a recursion relation
\begin{align}
x_i & = P_i + B_{i, \emptyset, \emptyset}, \notag \\
x_i P_u & = P_{(i, u)} + \sum_{w=1}^d B_{i, w, u} P_{w} + \delta_{i, u} C_u, \label{Recursion} \\
x_i P_{\vec{u}} & = P_{(i, \vec{u})} + \sum_{\abs{\vec{w}} = \abs{\vec{u}}} B_{i, \vec{w}, \vec{u}} P_{\vec{w}} + \delta_{i, u(1)} C_{\vec{u}} P_{(u(2), u(3), \ldots, u(k))}, \notag
\end{align}
with $C_{\vec{u}} \geq 0$ and, denoting $\vec{s}_j = (s(j), \ldots, s(k))$,
\[
B_{i, \vec{s}, \vec{u}} \prod_{j=1}^k C_{\vec{s}_j} = B_{i, \vec{u}, \vec{s}} \prod_{j=1}^k C_{\vec{u}_j}.
\]
\item
The state $\phi$ has a MOPS $\set{P_{\vec{u}}}$ and for any $\vec{u} \neq \vec{w}$, $\abs{\vec{u}} = \abs{\vec{w}} = n$,
\begin{equation}
\label{Relation-0}
\ip{x_{\vec{u}}}{{x_{\vec{w}}}}
= \sum_{\abs{\vec{v}} < n, \norm{P_{\vec{v}}} \neq 0}
\frac{\ip{x_{\vec{u}}}{P_{\vec{v}}} \ip{P_{\vec{v}}}{x_{\vec{w}}}}{\ip{P_{\vec{v}}}{P_{\vec{v}}}},
\end{equation}
so that the even non-symmetric moments of $\phi$ are determined by the rest of its moments.

\item
The state $\phi$ has a Fock space representation $\phi_{\mc{C}, \set{\mc{T}_i}}$ as in Definition~\ref{Defn:Fock-state}.
\end{enumerate}

\br
If $\phi$ is faithful, the following two conditions are equivalent to the preceding ones:
\begin{enumerate}
\item[(f)]
The polynomials $\set{\frac{1}{\mk{h}_{\abs{\vec{u}}}} \det M_{\vec{u}}(\mb{x})}$ (from Notation~\ref{Notation:Hankel}) are orthogonal with respect to $\phi$.
\item[(g)]
For any $\vec{u} \neq \vec{w}$, $\abs{\vec{u}} = \abs{\vec{w}} = n$,
\begin{equation}
\label{Relation-1}
\ip{x_{\vec{u}}}{{x_{\vec{w}}}}
= \sum_{\abs{\vec{v}} < n, h_{\vec{v}} \neq 0} \frac{h_{\vec{u}, \vec{v}} h_{\vec{w}, \vec{v}}}{h_{\vec{v}} \mk{h}_{\abs{\vec{v}}}}.
\end{equation}
\end{enumerate}
\end{Thm}

\begin{proof}
(a) $\Rightarrow$ (b).
Suppose that $\phi$ has a MOPS $\set{Q_{\vec{u}}}$. Proceed by induction on $n = \abs{\vec{u}}$. $P_{\emptyset} = Q_{\emptyset} = 1$. Suppose that $\ip{P_{\vec{w}}}{P_{\vec{v}}} = 0$ for all $\abs{\vec{w}}, \abs{\vec{v}} < n$. From relation~\eqref{Gram-Schmidt}, it follows that $\ip{P_{\vec{u}}}{P_{\vec{v}}} = 0$ for $\abs{\vec{v}} < \abs{\vec{u}}$. Using one direction of Lemma~\ref{Lemma:L2-equal}, we conclude that $\norm{P_{\vec{u}} - Q_{\vec{u}}} = 0$. Using the other direction of Lemma~\ref{Lemma:L2-equal}, it follows that for any $\vec{u} \neq \vec{v}$ with $\abs{\vec{u}}, \abs{\vec{v}} \leq n$, $\ip{P_{\vec{u}}}{P_{\vec{v}}} = 0$, and the induction hypothesis is satisfied. Therefore $\set{P_{\vec{u}}}$ are a MOPS.

\br
(b) $\Rightarrow$ (d).
By Lemma~\ref{Lemma:L2-equal}, all MOPS for $\phi$ have the same inner products, so we might as well use the MOPS from part (b) which satisfy the Gram-Schmidt recursion~\eqref{Gram-Schmidt}. Then for $\vec{u} \neq \vec{w}$, $\abs{\vec{u}} = \abs{\vec{w}} = n$,
\begin{equation}
\label{Relation}
\begin{split}
\ip{x_{\vec{u}}}{x_{\vec{w}}}
& = \ip{P_{\vec{u}}}{P_{\vec{w}}}
+ \sum_{\vec{v}:\ \abs{\vec{v}} < \abs{\vec{w}},\ \norm{P_{\vec{v}}} \neq 0} \frac{\ip{x_{\vec{w}}}{P_{\vec{v}}} \ip{P_{\vec{u}}}{P_{\vec{v}}}}{\norm{P_{\vec{v}}}^2}
+ \sum_{\vec{v}:\ \abs{\vec{v}} < \abs{\vec{u}},\ \norm{P_{\vec{v}}} \neq 0} \frac{\ip{x_{\vec{u}}}{P_{\vec{v}}} \ip{P_{\vec{v}}}{P_{\vec{w}}}}{\norm{P_{\vec{v}}}^2} \\
&\quad + \sum_{\vec{v}, \vec{s} :\ \abs{\vec{v}}, \abs{\vec{s}} < \abs{\vec{w}},\ \norm{P_{\vec{v}}} \neq 0 \neq \norm{P_{\vec{s}}}} \frac{\ip{x_{\vec{u}}}{P_{\vec{v}}}}{\norm{P_{\vec{v}}}^2} \frac{\ip{x_{\vec{w}}}{P_{\vec{s}}}}{\norm{P_{\vec{s}}}^2} \ip{P_{\vec{v}}}{P_{\vec{s}}} \\
& = \sum_{\abs{\vec{v}} < n, \norm{P_{\vec{v}}} \neq 0} \frac{\ip{x_{\vec{u}}}{P_{\vec{v}}} \ip{P_{\vec{v}}}{x_{\vec{w}}}}{\ip{P_{\vec{v}}}{P_{\vec{v}}}}.
\end{split}
\end{equation}

\br
(d) $\Rightarrow$ (a).
Obvious.

\br
(a) $\Rightarrow$ (e).
Let $\phi$ be a state with MOPS $\set{P_{\vec{u}}}$. For each $k$, define a diagonal matrix $\mc{C}^{(k)}$ recursively via
\begin{equation}
\label{C1}
\begin{split}
\state{P_{\vec{u}}^\ast(\mb{x}) P_{\vec{u}}(\mb{x})}
& = \ip{e_{u(1)} \otimes \mc{K}_{\mc{C}} (e_{u(2)} \otimes \ldots \otimes e_{u(k)})}{\mc{C}^{(k)} (e_{u(1)} \otimes \ldots \otimes e_{u(k)})} \\
& = \norm{e_{u(1)} \otimes \ldots \otimes e_{u(k)}}_{\mc{C}}.
\end{split}
\end{equation}
This may not be well-defined if $\mc{K}_{\mc{C}}(e_{u(2)} \otimes \ldots \otimes e_{u(k)}) = 0$. However, in that case $\norm{P_{(u(2), \ldots u(k))}} = 0$, which by Lemma~\ref{Lemma:Ideal} implies that $\state{P_{\vec{u}}^\ast(\mb{x}) P_{\vec{u}}(\mb{x})} = 0$, so one can take the corresponding entry of $\mc{C}^{(k)}$ to be anything, for example zero. Similarly, for each $k$ and for $\abs{\vec{u}} = \abs{\vec{v}}=k$, define $\mc{T}_i^{(k)}$ via
\begin{equation}
\label{T1}
\begin{split}
\state{P_{\vec{u}}^\ast(\mb{x}) x_i P_{\vec{v}}(\mb{x})}
& = \ip{\mc{K}_{\mc{C}} (e_{u(1)} \otimes \ldots \otimes e_{u(k)})}{\mc{T}_i (e_{v(1)} \otimes \ldots \otimes e_{v(k)})} \\
& = \ip{e_{u(1)} \otimes \ldots \otimes e_{u(k)}}{\mc{T}_i (e_{v(1)} \otimes \ldots \otimes e_{v(k)})}_{\mc{C}}.
\end{split}
\end{equation}
Again, this may not be well-defined if $\norm{e_{u(1)} \otimes \ldots \otimes e_{u(k)}}_{\mc{C}} = 0$ or equivalently if $\norm{P_{\vec{u}}(\mb{x})} = 0$. But in that case $\state{P_{\vec{u}}^\ast(\mb{x}) x_i P_{\vec{v}}(\mb{x})} = \ip{P_{\vec{u}}}{x_i P_{\vec{v}}} = 0$, and one can take the corresponding entry of $\mc{T}_i^{(k)}$ to be anything, for example zero.

\br
By construction, each $\mc{C}^{(k)}$ is diagonal and non-negative. Also,
\[
\state{P_{\vec{u}}^\ast(\mb{x}) x_i P_{\vec{v}}(\mb{x})}
= \state{\bigl(x_i P_{\vec{u}}(\mb{x})\bigr)^\ast P_{\vec{v}}(\mb{x})}
= \state{P_{\vec{v}}^\ast (\mb{x}) x_i P_{\vec{u}}(\mb{x})}
\]
so $(\mc{K}_{\mc{C}} \mc{T}_i)^t = \mc{T}_i^t \mc{K}_{\mc{C}} = \mc{K}_{\mc{C}} \mc{T}_i$.

\br
Let $\set{Q_{\vec{u}}}$ be a MOPS for $\phi_{\mc{C}, \set{\mc{T}_i}}$ constructed in Lemma~\ref{Lemma:Fock-polynomials}. Then for any $\abs{\vec{u}} = \abs{\vec{v}}$, the equations~\eqref{T1} and~\eqref{C1} translate to
\begin{equation}
\label{C2}
\state{P_{\vec{u}}^\ast(\mb{x}) P_{\vec{v}}(\mb{x})}
= \phi_{\mc{C}, \set{\mc{T}_i}}\left[Q_{\vec{u}}^\ast(\mb{x}) Q_{\vec{v}}(\mb{x})\right]
\end{equation}
and
\begin{equation}
\label{T2}
\begin{split}
\state{P_{\vec{u}}^\ast(\mb{x}) x_i P_{\vec{v}}(\mb{x})}
= \phi_{\mc{C}, \set{\mc{T}_i}}\left[Q_{\vec{u}}^\ast(\mb{x}) x_i Q_{\vec{v}}(\mb{x})\right].
\end{split}
\end{equation}
We may assume that both $\set{P_{\vec{u}}}$ and $\set{Q_{\vec{u}}}$ satisfy, for their respective states, the Gram-Schmidt recursions~\eqref{Gram-Schmidt}. Now we prove, by induction on the degree, that one can take $Q_{\vec{u}} = P_{\vec{u}}$, which implies the equality of the states $\phi = \phi_{\mc{C}, \set{\mc{T}_i}}$. Suppose that $Q_{\vec{v}} = P_{\vec{v}}$ for $\abs{\vec{u}} \leq n$. Then it follows from equations~\eqref{C2} and~\eqref{T2} that the two states coincide on polynomials of degree at most $2n+1$. But in that case, the Gram-Schmidt recursions for $P_{\vec{u}}$ and $Q_{\vec{u}}$, for $\abs{\vec{u}} = n+1$, are identical.

\br
(e) $\Rightarrow$ (c).
Let $\phi = \phi_{\mc{C}, \set{\mc{T}_i}}$. Let $\set{B_{i, \vec{w}, \vec{u}}, C_{\vec{u}}}$ be the matrix elements of $\set{\mc{T}_i^{(k)}, \mc{C}^{(k)}}$, so that
\begin{equation*}
\mc{T}_i(e_{u(1)} \otimes \ldots \otimes e_{u(k)}) = \sum_{\abs{\vec{w}} = k} B_{i, \vec{w}, \vec{u}} e_{w(1)} \otimes \ldots \otimes e_{w(k)}
\end{equation*}
and
\begin{equation*}
\mc{C}(e_{u(1)} \otimes \ldots \otimes e_{u(k)}) = C_{\vec{u}} e_{u(1)} \otimes \ldots \otimes e_{u(k)}.
\end{equation*}
Note that the conditions on the coefficients in part (c) correspond exactly to the conditions on operators in Section~\ref{Subsec:General-Fock}. Define the polynomials $\set{P_{\vec{u}}}$ using the recursion~\eqref{Recursion}. We show by induction that relation~\eqref{Wick-representation} holds. Indeed,
\[
\begin{split}
P_{(i, \vec{u})}(\mb{X}) \Omega
& = \mc{X}_i P_{\vec{u}}(\mb{X}) \Omega - \sum_{\abs{\vec{w}} = \abs{\vec{u}}} B_{i, \vec{w}, \vec{u}} P_{\vec{w}}(\mb{X}) \Omega - \delta_{i, u(1)} C_{\vec{u}} P_{(u(2), u(3), \ldots, u(k))}(\mb{X}) \Omega \\
& = \mc{X}_i e_{u(1)} \otimes \ldots \otimes e_{u(k)} - \sum_{\abs{\vec{w}} = \abs{\vec{u}}} B_{i, \vec{w}, \vec{u}} e_{w(1)} \otimes \ldots \otimes e_{w(k)} - \delta_{i, u(1)} C_{\vec{u}} e_{u(2)} \otimes \ldots \otimes e_{u(k)} \\
& = \mc{X}_i (e_{u(1)} \otimes \ldots \otimes e_{u(k)}) - \mc{T}_i(e_{u(1)} \otimes \ldots \otimes e_{u(k)}) - a_i^- \mc{C} (e_{u(1)} \otimes \ldots \otimes e_{u(k)}) \\
& = a_i^+ (e_{u(1)} \otimes \ldots \otimes e_{u(k)})
= e_i \otimes e_{u(1)} \otimes \ldots \otimes e_{u(k)}.
\end{split}
\]
It follows that for all $\vec{u} \neq \emptyset$, $\state{P_{\vec{u}}(\mb{x})} = 0$.

\br
(c) $\Rightarrow$ (a).
For this direction only, the arguments of Proposition 3 of \cite{AnsMulti-Sheffer} apply, and it follows that $\set{P_{\vec{u}}}$ are orthogonal. Note that in that proposition $\phi$ was assumed to have zero means and identity covariance; however, it is easy to modify its proof to get a state even without those assumptions.

\br
Now assume that a state $\phi$ is faithful. Then for every $n$, $\mk{h}_n \neq 0$. Indeed, if some $\mk{h}_n = 0$, then for the corresponding matrix, a linear combination of some of its rows is zero. This is equivalent to saying that for some polynomial $P$ of degree less than $n$, $\ip{P}{x_{\vec{v}}} = 0$ for all $\abs{\vec{v}} < n$. But in that case, $\norm{P} = 0$, which contradicts the faithfulness assumption.

\br
$\mk{h}_{\abs{\vec{u}}}$ is the leading coefficient of $\det M_{\vec{u}}(\mb{x})$, and since it is non-zero, $\set{\frac{1}{\mk{h}_{\abs{\vec{u}}}} \det M_{\vec{u}}}$ are a monic polynomial family. Moreover, for $\abs{\vec{v}} < \abs{\vec{u}}$, $\ip{x_{\vec{v}}}{\det M_{\vec{u}}}$ equals $h_{\vec{v}, \vec{u}}$, the determinant of the matrix $A_{\vec{u}}$ with the $\vec{u}$'th row replaced so that the $(\vec{u}, \vec{w})$ entry is $\ip{x_{\vec{v}}}{x_{\vec{w}}}$. Such a row is identical to the $\vec{v}$'th row, and so the determinant is zero. Thus these polynomials are pseudo-orthogonal.

\br
(a) $\Leftrightarrow$ (f).
One implication follows from Lemma~\ref{Lemma:L2-equal}. The other one is trivial.

\br
(f) $\Leftrightarrow$ (g).
Given $\phi$, let $\set{P_{\vec{u}}}$ be a pseudo-orthogonal family defined via the Gram-Schmidt recursion~\eqref{Gram-Schmidt}. From the calculation~\eqref{Relation}, it follows that
\[
\ip{x_{\vec{u}}}{x_{\vec{w}}}
= \ip{P_{\vec{u}}}{P_{\vec{w}}}
+ \sum_{\abs{\vec{v}} < n, \norm{P_{\vec{v}}} \neq 0} \frac{\ip{x_{\vec{u}}}{P_{\vec{v}}} \ip{P_{\vec{v}}}{x_{\vec{w}}}}{\ip{P_{\vec{v}}}{P_{\vec{v}}}}.
\]
Thus for pseudo-orthogonal polynomials, condition~\eqref{Relation-0} is equivalent to orthogonality. It remains to note that
\[
\ip{x_{\vec{u}}}{\det M_{\vec{v}}} = h_{\vec{u}, \vec{v}}
\]
and in particular
\[
\ip{x_{\vec{v}}}{\det M_{\vec{v}}} = h_{\vec{v}}
\]
so that
\[
\norm{\frac{1}{\mk{h}_{\abs{\vec{v}}}} \det M_{\vec{v}}}^2 = \frac{h_{\vec{v}}}{\mk{h}_{\abs{\vec{v}}}},
\]
and condition~\eqref{Relation-1} is exactly the condition~\eqref{Relation-0} for the polynomials $\set{\frac{1}{\mk{h}_{\abs{\vec{u}}}} \det M_{\vec{u}}}$.
\end{proof}

\begin{Ex}
For degree one, parts (d, g) of the theorem say that for $i \neq j$,
\[
\state{x_i x_i} = \state{x_i} \state{x_j},
\]
so that the variables are uncorrelated. For degree two, they say that for $(i,j) \neq (t,s)$, and assuming for simplicity that $\state{x_i} = 0$ for all $i$,
\[
\state{x_i x_j x_s x_t} = \state{x_i x_j} \state{x_s x_t} + \sum_{k=1}^d \frac{\state{x_i x_j x_k} \state{x_k x_s x_t}}{\state{x_k^2}}.
\]
\end{Ex}


\providecommand{\bysame}{\leavevmode\hbox to3em{\hrulefill}\thinspace}
\providecommand{\MR}{\relax\ifhmode\unskip\space\fi MR }
\providecommand{\MRhref}[2]{%
  \href{http://www.ams.org/mathscinet-getitem?mr=#1}{#2}
}
\providecommand{\href}[2]{#2}

\end{document}